\documentclass[leqno,draft]{article}

\newtheorem{theorem}{Theorem}
\newtheorem{lemma}[theorem]{Lemma}
\newtheorem{proposition}[theorem]{Proposition}
\newtheorem{definition}[theorem]{Definition}
\newtheorem{corollary}[theorem]{Corollary}

\newcommand{\begintheorem}{\addtocounter{equation}{1}\begin{theorem}}
\newcommand{\beginlemma}{\addtocounter{equation}{1}\begin{lemma}}
\newcommand{\beginproposition}{\addtocounter{equation}{1}\begin{proposition}}
\newcommand{\begindefinition}{\addtocounter{equation}{1}\begin{definition}}
\newcommand{\begincorollary}{\addtocounter{equation}{1}\begin{corollary}}

\begin{document}

\title{Another introduction to the geometry \\ of metric spaces}

\author{Stephen Semmes \\
        Rice University}

\date{}

\maketitle


\renewcommand{\thefootnote}{}   

\footnotetext{The present notes are somewhat complementary to
\cite{s2}, and the two can be read in either order.}

\begin{abstract}
Here Lipschitz conditions are used as a primary tool, for studying
curves in metric spaces in particular.
\end{abstract}

\tableofcontents

\section{Basic notions}
\label{basic notions}
\setcounter{equation}{0}

        A \emph{metric space} is a nonempty set $M$ equipped with a
distance function $d(x, y)$ defined for $x, y \in M$ such that $d(x,
y)$ is a nonnegative real number which is equal to $0$ exactly when $x
= y$,
\begin{equation}
        d(y, x) = d(x, y)
\end{equation}
for every $x, y \in M$, and
\begin{equation}
        d(x, z) \le d(x, y) + d(y, z)
\end{equation}
for every $x, y, z \in M$.  This last condition is known as the
\emph{triangle inequality}.

        As usual, ${\bf R}$ denotes the real line, and the
\emph{absolute value} $|r|$ of $r \in {\bf R}$ is defined to be $r$
when $r \ge 0$ and $-r$ when $r \le 0$.  It is easy to check that
\begin{equation}
\label{|r + t| le |r| + |t|}
        |r + t| \le |r| + |t|
\end{equation}
for every $r, t \in {\bf R}$, which implies that $|r - t|$ is a metric
on ${\bf R}$.

        Let $(M, d(x, y))$ be a metric space.  For each $x \in M$ and
$r > 0$, the \emph{open ball} with center $x$ and radius $r$ is defined by
\begin{equation}
        B(x, r) = \{y \in M : d(x, y) < r\},
\end{equation}
and the corresponding \emph{closed ball} is defined to be
\begin{equation}
        \overline{B}(x, r) = \{y \in M : d(x, y) \le r\}.
\end{equation}

        A set $E \subseteq M$ is said to be \emph{bounded} if $E$ is
contained in a ball in $M$.  For any $p, q \in M$ and $r > 0$,
\begin{equation}
        B(p, r) \subseteq B(q, r + d(p, q))
\end{equation}
and
\begin{equation}
        \overline{B}(p, r) \subseteq \overline{B}(q, r + d(p, q)),
\end{equation}
by the triangle inequality.  It follows that a bounded set $E
\subseteq M$ is contained in a ball centered at any point in $M$.

\section{Norms on ${\bf R}^n$}
\label{norms}
\setcounter{equation}{0}

        For each positive integer $n$, ${\bf R}^n$ is the space of
$n$-tuples $x = (x_1, \ldots, x_n)$ of real numbers, i.e., $x_1,
\ldots, x_n \in {\bf R}$.  This is a vector space with respect to
coordinatewise addition and scalar multiplication by real numbers.

        Suppose that $N(x)$ is a function defined on ${\bf R}^n$ with
values in the nonnegative real numbers.  We say that $N(x)$ is a
\emph{norm} on ${\bf R}^n$ if $N(x) = 0$ exactly when $x = 0$,
\begin{equation}
\label{N(r x) = |r| N(x)}
        N(r \, x) = |r| \, N(x)
\end{equation}
for every $r \in {\bf R}$ and $x \in {\bf R}^n$, and
\begin{equation}
\label{N(x + y) le N(x) + N(y)}
        N(x + y) \le N(x) + N(y)
\end{equation}
for every $x, y \in {\bf R}^n$.  If $N(x)$ is a norm on ${\bf R}^n$,
then
\begin{equation}
\label{d_N(x, y) = N(x - y)}
        d_N(x, y) = N(x - y)
\end{equation}
is a metric on ${\bf R}^n$.

        The absolute value function is a norm on the real line, and
any norm on ${\bf R}$ can be expressed as $a \, |x|$ for some $a > 0$.
The standard Euclidean norm on ${\bf R}^n$ is defined by
\begin{equation}
        |x| = \Big(\sum_{j = 1}^n x_j^2\Big)^{1/2},
\end{equation}
and determines the standard metric on ${\bf R}^n$.  It is not
completely obvious that this satisfies the triangle inequality, and
one way to show this will be mentioned in the next section.  It is
much easier to check directly that
\begin{equation}
        \|x\|_1 = \sum_{j = 1}^n |x_j|
\end{equation}
and
\begin{equation}
        \|x\|_\infty = \max(|x_1|, \ldots, |x_n|)
\end{equation}
are norms on ${\bf R}^n$.  Note that the standard norm may also be
denoted $\|x\|_2$.

         If $N(x)$ is any norm on ${\bf R}^n$, then
\begin{equation}
         N(x) \le \max(N(e_1), \ldots, N(e_n)) \, \|x\|_1,
\end{equation}
where $e_1, \ldots, e_n$ are the standard basis vectors in ${\bf
R}^n$, which is to say that the $j$th coordinate of $e_\ell$ is equal
to $1$ when $j = \ell$ and to $0$ otherwise.  Indeed,
\begin{equation}
         x = \sum_{j = 1}^n x_j \, e_j,
\end{equation}
and therefore
\begin{equation}
         N(x) \le \sum_{j = 1}^n N(e_j) \, |x_j|.
\end{equation}
In particular,
\begin{equation}
         \|x\|_\infty \le \|x\|_2 \le \|x\|_1
\end{equation}
for every $x \in {\bf R}^n$, since the first inequality holds by
inspection.  One can also get the second inequality by observing that
\begin{equation}
         \sum_{j = 1}^n x_j^2 \le \|x\|_1 \, \|x\|_\infty,
\end{equation}
since $\|x\|_\infty \le \|x\|_1$ also holds by inspection.  In the
other direction, it is easy to see that
\begin{equation}
        \|x\|_1 \le n \, \|x\|_\infty
\end{equation}
and
\begin{equation}
         \|x\|_2 \le \sqrt{n} \, \|x\|_\infty,
\end{equation}
and one can use the convexity of $\phi(t) = t^2$ on the real line to show that
\begin{equation}
         \|x\|_1 \le \sqrt{n} \, \|x\|_2.
\end{equation}

\section{Convex sets in ${\bf R}^n$}
\label{convex sets in R^n}
\setcounter{equation}{0}

        A set $E \subseteq {\bf R}^n$ is said to be \emph{convex} if
for every $x, y \in E$ and every real number $t$, $0 \le t \le 1$,
\begin{equation}
        t \, x + (1 - t) \, y \in E.
\end{equation}
For example, open and closed balls with respect to the metric
associated to a norm on ${\bf R}^n$ are convex.

        Conversely, suppose that $N(x)$ is a nonnegative real-valued
function on ${\bf R}^n$ that satisfies $N(x) > 0$ when $x \ne 0$ and
the homogeneity condition (\ref{N(r x) = |r| N(x)}).  If
\begin{equation}
\label{B_N = {x in {bf R}^n : N(x) le 1}}
        B_N = \{x \in {\bf R}^n : N(x) \le 1\}
\end{equation}
is convex, then $N$ satisfies the triangle inequality (\ref{N(x + y)
le N(x) + N(y)}) and hence is a norm.  To see this, let $x, y \in {\bf
R}^n$ be given, and let us check (\ref{N(x + y) le N(x) + N(y)}).  The
inequality is trivial when $x = 0$ or $y = 0$, and so we may suppose
that $x, y \ne 0$.  Put
\begin{equation}
        x' = \frac{x}{N(x)}, \quad y' = \frac{y}{N(y)},
\end{equation}
which automatically satisfy
\begin{equation}
        N(x') = N(y') = 1.
\end{equation}
For $0 \le t \le 1$, convexity of $B_N$ implies that
\begin{equation}
\label{N(t x' + (1 - t) y') le 1}
        N(t \, x' + (1 - t) \, y') \le 1.
\end{equation}
If
\begin{equation}
        t = \frac{N(x)}{N(x) + N(y)},
\end{equation}
then
\begin{equation}
        1 - t = \frac{N(y)}{N(x) + N(y)}
\end{equation}
and
\begin{equation}
        t \, x' + (1 - t) \, y' = \frac{x + y}{N(x) + N(y)},
\end{equation}
which means that (\ref{N(x + y) le N(x) + N(y)}) follows from
(\ref{N(t x' + (1 - t) y') le 1}).

        One can use the convexity of the function $\phi(r) = r^2$ on
the real line to show directly that the closed unit ball with respect
to the standard Euclidean norm is a convex set in ${\bf R}^n$, and
hence that the Euclidean norm satisfies the triangle inequality and is
therefore a norm.  For each real number $p$, $1 \le p < \infty$, put
\begin{equation}
        \|x\|_p = \Big(\sum_{j = 1}^n |x_j|^p \Big)^{1/p}.
\end{equation}
One can use the convexity of the function $\phi_p(r) = |r|^p$ on the
real line to show that the closed unit ball associated to $\|x\|_p$ is
a convex set in ${\bf R}^n$, and therefore that $\|x\|_p$ is a norm on
${\bf R}^n$.

        By inspection,
\begin{equation}
        \|x\|_\infty \le \|x\|_p
\end{equation}
for every $x \in {\bf R}^n$ and $1 \le p < \infty$.  If $1 \le p < q <
\infty$, then
\begin{equation}
        \sum_{j = 1}^n |x_j|^q
          \le \Big(\sum_{j = 1}^n |x_j|^p \Big) \, \|x\|_\infty^{q - p},
\end{equation}
which implies that
\begin{equation}
        \|x\|_q \le \|x\|_p^{p/q} \, \|x\|_\infty^{1 - (p/q)}
\end{equation}
and thus
\begin{equation}
        \|x\|_q \le \|x\|_p
\end{equation}
for every $x \in {\bf R}^n$.

\section{Lipschitz conditions, 1}
\label{lipschitz conditions, 1}
\setcounter{equation}{0}

        Let $(M_1, d_1(x, y))$, $(M_2, d_2(u, v))$ be metric spaces.  A
mapping $f : M_1 \to M_2$ is said to be \emph{Lipschitz} with constant
$C \ge 0$ or $C$-Lipschitz if
\begin{equation}
        d_2(f(x), f(y)) \le C \, d_1(x, y)
\end{equation}
for every $x, y \in M_1$.  More precisely, $f$ is Lipschitz of order
$1$ if this holds for some $C \ge 0$, and we shall discuss Lipschitz
conditions of any order $a > 0$ a bit later.  Note that $f$ is
Lipschitz with $C = 0$ if and only if $f$ is constant, and that
Lipschitz mappings are automatically uniformly continuous.

        If $M_2$ is the real line with the standard metric,
then the preceding Lipschitz condition is equivalent to
\begin{equation}
        f(x) \le f(y) + C \, d_1(x, y).
\end{equation}
This follows by interchanging the order of $x$ and $y$.  In particular,
\begin{equation}
        f_p(x) = d_1(x, p)
\end{equation}
is Lipschitz with $C = 1$ for every $p \in M_1$, by the triangle
inequality.  If $f$, $\widetilde{f}$ are real-valued Lipschitz
functions on $M_1$ with constants $C$, $\widetilde{C}$, respectively,
then $f + \widetilde{f}$ is Lipschitz with constant $C +
\widetilde{C}$.  Moreover, $a \, f$ is Lipschitz with constant $|a| \,
C$ for every $a \in {\bf R}$.  The product of bounded real-valued
Lipschitz functions is also Lipschitz.  If $(M_3, d_3(w, z))$ is
another metric space, and $f_1 : M_1 \to M_2$ and $f_2 : M_2 \to M_3$
are Lipschitz mappings with constants $C_1$, $C_2$, respectively, then
the composition $f_2 \circ f_1 : M_1 \to M_3$ defined by
\begin{equation}
        (f_2 \circ f_1)(x) = f_2(f_1(x))
\end{equation}
is Lipschitz with constant $C_1 \, C_2$.

        For any mapping $f : M_1 \to M_2$ and set $A \subseteq M_1$,
\begin{equation}
        f(A) = \{f(x) : x \in A\} \subseteq M_2.
\end{equation}
Let $B_1(x, r)$ and $B_2(p, t)$ be the open balls in $M_1$, $M_2$ with
centers $x \in M_1$, $p \in M_2$ and radii $r, t > 0$, respectively.
It is easy to see that $f : M_1 \to M_2$ is Lipschitz with constant $C
> 0$ if and only if
\begin{equation}
        f(B_1(x, r)) \subseteq B_2(f(x), C \, r)
\end{equation}
for every $x \in M_1$ and $r > 0$.  This is also equivalent to the
analogous condition
\begin{equation}
        f(\overline{B}_1(x, r)) \subseteq \overline{B}_2(f(x), C \, r)
\end{equation}
for closed balls.  In particular, if $A$ is a bounded set in $M_1$,
then $f(A)$ is bounded in $M_2$.

        Suppose that $f$ is a real-valued function on the real line,
equipped with the standard metric.  If $f$ is differentiable at a
point $x \in {\bf R}$, and $f$ is $C$-Lipschitz for some $C \ge 0$,
then the derivative $f'(x)$ of $f$ at $x$ satisfies
\begin{equation}
\label{|f'(x)| le C}
        |f'(x)| \le C.
\end{equation}
This follows from the definition of the derivative.  Conversely, if
$f$ is differentiable and satisfies (\ref{|f'(x)| le C}) at every
point in ${\bf R}$, then $f$ is $C$-Lipschitz, by the mean value
theorem.  Note that $f(x) = |x|$ is $1$-Lipschitz on ${\bf R}$ and not
differentiable at $x = 0$.

\section{Lipschitz curves}
\label{lipschitz curves}
\setcounter{equation}{0}

        Let $(M, d(x, y))$ be a metric space, and suppose that $a$,
$b$ are real numbers with $a \le b$.  As usual, the closed interval
$[a, b]$ in the real line consists of the $r \in {\bf R}$ such that $a
\le r \le b$.  Suppose also that $p : [a, b] \to M$ is Lipschitz with
constant $k$ for some $k \ge 0$.  If $\{t_j\}_{j = 0}^n$ is a finite
sequence of real numbers such that
\begin{equation}
        a = t_0 < t_1 < \cdots < t_n = b,
\end{equation}
then
\begin{equation}
        \sum_{j = 1}^n d(p(t_j), p(t_{j - 1}))
          \le k \, \sum_{j = 1}^n (t_j - t_{j - 1})
           = k \, (b - a).
\end{equation}
This is often described by saying that the length of the curve
determined by $p(t)$, $a \le t \le b$, has length $\le k \, (b - a)$.

        Of course, one can use translations on the real line to shift
the interval on which a path is defined without changing the Lipschitz
constant.  One can use affine mappings on ${\bf R}$ to change the
length of the interval on which a path is defined, with a
corresponding change in the Lipschitz constant.  The product of the
Lipschitz constant and the length of the interval would remain the
same.

        If $c \in {\bf R}$, $c \ge b$, $q : [b, c] \to M$ is
$k$-Lipschitz, and $p(b) = q(b)$, then the mapping from $[a, c]$ into
$M$ defined by combining $p$ and $q$ is $k$-Lipschitz too.  This is
easy to verify, directly from the definitions.  If the Lipschitz
constants for $p$ and $q$ are different, then it may be preferable to
rescale the intervals so that the Lipschitz constants are the same.
If $p$ is constant on an interval $[a_1, b_1] \subseteq [a, b]$, then
one can remove $(a_1, b_1)$ from $[a, b]$ and combine the remaining
pieces to get a curve with the same Lipschitz constant on a smaller
interval.

\section{Minimality}
\label{minimality}
\setcounter{equation}{0}

        Let $(M, d(x, y))$ be a metric space in which closed and
bounded sets are compact.  Suppose that $x, y \in M$ can be connected
by a Lipschitz curve in $M$.  This means that there is a Lipschitz
mapping $p : [0, 1] \to M$ such that $p(0) = x$ and $p(1) = y$.  Using
the Arzela-Ascoli theorem, one can show that there is such a path
whose Lipschitz constant is as small as possible.  For suppose that
$p_1, p_2, \ldots$ is a sequence of Lipschitz mappings from $[0, 1]$
into $M$ whose Lipschitz constants $k_1, k_2, \ldots$, respectively,
converge to the infimum $k$ of the possible Lipschitz constants.  By
passing to a subsequence, we may suppose that the sequence of mappings
converges uniformly on $[0, 1]$.  The limiting mapping sends $0$ to
$x$ and $1$ to $y$, and it is easy to check that it is Lipschitz with
constant $k$.  Note that $k \ge d(x, y)$.

\section{Affine paths in ${\bf R}^n$}
\label{affine paths in R^n}
\setcounter{equation}{0}

        Fix a positive integer $n$, and consider an affine mapping $p
: {\bf R} \to {\bf R}^n$, given by $p(r) = u + r \, v$ for some $u, v
\in {\bf R}^n$.  If $N$ is any norm on ${\bf R}^n$, then $p$ is
Lipschitz with constant $N(v)$ with respect to the standard metric on
${\bf R}$ and the metric associated to $N$ on ${\bf R}^n$, since
\begin{equation}
        N(p(r) - p(t)) = |r - t| \, N(v)
\end{equation}
for every $r, t \in {\bf R}^n$.  For any $a, b \in {\bf R}$ with $a
\le b$, the restriction of $p$ to $[a, b]$ is a Lipschitz curve
connecting $p(a)$ to $p(b)$ with constant $N(v)$, and $N(v)$ is the
smallest possible Lipschitz constant for such a curve on $[a, b]$.

        A norm $N$ on ${\bf R}^n$ is said to be \emph{strictly convex}
if the corresponding closed unit ball $B_N$ as in (\ref{B_N = {x in
{bf R}^n : N(x) le 1}}) is strictly convex.  This means that for every
$x, y \in {\bf R}^n$ with $N(x) = N(y) = 1$ and $x \ne y$, we have
that
\begin{equation}
        N(t \, x + (1 - t) \, y) < 1
\end{equation}
when $0 < t < 1$.  Equivalently, if $w, z \in {\bf R}^n$ and
\begin{equation}
        N(w + z) = N(w) + N(z),
\end{equation}
then either $w = 0$, $z = 0$, or $z = r \, w$ for some $r > 0$.  This
follows from an argument like the one used in Section \ref{convex sets
in R^n} to show that convexity of $B_N$ implies the triangle
inequality for $N$ when $N$ is homogeneous.

        One can show that the standard Euclidean norm on ${\bf R}^n$
is strictly convex, using strict convexity of the function $\phi(r) =
r^2$ on the real line.  Similarly, $\|x\|_p$ is strictly convex on
${\bf R}^n$ when $1 < p < \infty$, as a consequence of the strict
convexity of $\phi_p(r) = |r|^p$.  In particular, the absolute value
is strictly convex as a norm on ${\bf R}$, if not in the ordinary
sense for arbitrary functions, because equality holds in (\ref{|r + t|
le |r| + |t|}) only when $r, t \ge 0$ or $r, t \le 0$.  However,
$\|x\|_1$ and $\|x\|_\infty$ are not strictly convex norms on ${\bf
R}^n$ when $n \ge 2$.

        Suppose that $N$ is a strictly convex norm on ${\bf R}^n$.
If $x, y, z \in {\bf R}^n$ satisfy
\begin{equation}
        d_N(x, z) = d_N(x, y) + d_N(y, z),
\end{equation}
where $d_N$ is as defined in (\ref{d_N(x, y) = N(x - y)}), then
\begin{equation}
        y = r \, x + (1 - r) \, z
\end{equation}
for some $r \in [0, 1]$.  If $q : [a, b] \to {\bf R}^n$ is
$k$-Lipschitz and
\begin{equation}
        N(q(b) - q(a)) = k \, (b - a),
\end{equation}
then
\begin{equation}
        d_N(q(a), q(b)) = d_N(q(a), q(t)) + d_N(q(t), q(b))
\end{equation}
for every $t \in [a, b]$.  One can use this to show that $q(t)$ is
affine.  This does not work for the norms $\|x\|_1$, $\|x\|_\infty$ on
${\bf R}^n$ when $n \ge 2$.  For example, there is a $1$-Lipschitz
path from $[0, 2]$ into ${\bf R}^2$ equipped with the norm $\|x\|_1$
that connects $(0, 0)$ to $(1, 1)$ by following the horizontal segment
to $(1, 0)$ and then the vertical segment to $(1, 1)$.
If $\phi : [0, 1] \to {\bf R}$ is any $1$-Lipschitz function with
respect to the standard metric on the real line which satisfies
$\phi(0) = \phi(1) = 0$, then $\Phi(t) = (t, \phi(t))$ is a
$1$-Lipschitz mapping from $[0, 1]$ into ${\bf R}^2$ equipped with
the norm $\|x\|_\infty$ that connects $(0, 0)$ to $(1, 0)$.

\section{$C^1$ paths in ${\bf R}^n$}
\label{C^1 paths in R^n}
\setcounter{equation}{0}

       Let $N$ be a norm on ${\bf R}^n$.  As in Section \ref{lipschitz
conditions, 1}, the triangle inequality implies that $N$ is
$1$-Lipschitz with respect to the associated metric $d_N$.  As in
Section \ref{norms}, one can show that $N$ is less than or equal to a
constant multiple of the standard Euclidean norm on ${\bf R}^n$.  It
follows that $N$ is also a Lipschitz function with respect to the
standard metric on ${\bf R}^n$.

       Let $p : [a, b] \to {\bf R}^n$ be a continuously-differentiable
curve with derivative $p'(t)$.  This implies that $N(p'(t))$ is a
continuous function on $[a, b]$.  By the fundamental theorem of calculus,
\begin{equation}
       p(t) - p(r) = \int_r^t p'(u) \, du
\end{equation}
when $a \le r \le t \le b$.  Hence
\begin{equation}
       N(p(t) - p(r)) \le \int_r^t N(p'(u)) \, du,
\end{equation}
using an extension of the triangle inequality from sums to integrals.  If
\begin{equation}
\label{N(p'(u)) le k}
       N(p'(u)) \le k
\end{equation}
for every $u \in [a, b]$, then it follows that $p$ is $k$-Lipschitz
with respect to the metric associated to $N$ on ${\bf R}^n$.

        Alternatively, let $\epsilon > 0$ be given.  For each $r \in
[a, b]$,
\begin{equation}
        p(t) - p(r) - p'(r) \, (t - r)
\end{equation}
is $\epsilon$-Lipschitz as a function of $t$ on sufficiently small
neighborhoods of $r$ in $[a, b]$, since $p$ is continuously-differentiable.
Under the hypothesis (\ref{N(p'(u)) le k}), we get that $p$ is $(k +
\epsilon)$-Lipschitz with respect to $N$ on sufficiently small neighborhoods
of every point in $[a, b]$.  One can use this to show that $p$ is
$(k + \epsilon)$-Lipschitz on $[a, b]$, and therefore $k$-Lipschitz
because $\epsilon > 0$ is arbitrary.  Note that (\ref{N(p'(u)) le k})
holds when $p$ is $k$-Lipschitz with respect to $N$ on $[a, b]$.

        In order for the product of the Lipschitz constant and the
length of the parameter interval to be as small as possible, it would
be nice to have $N(p')$ constant on $[a, b]$.  As in the classical
situation, one can try to get this by reparameterizing $p$.  This is
easy to do when $p'(t) \ne 0$ for every $t \in [a, b]$.  Specifically,
\begin{equation}
        \phi(t) = \int_a^t N(p'(u)) \, du
\end{equation}
is a continuously-differentiable function on $[a, b]$ with
\begin{equation}
        \phi'(t) = N(p'(t)) > 0
\end{equation}
for each $t \in [a, b]$.  If $q = p \circ \phi^{-1}$, then
\begin{equation}
        N(q'(r)) = 1
\end{equation}
when $\phi(a) \le r \le \phi(b)$.

\section{Lipschitz conditions, 2}
\label{lipschitz conditions, 2}
\setcounter{equation}{0}

        Let $(M_1, d_1(x, y))$ and $(M_2, d_2(u, v))$ be metric
spaces.  A mapping $f : M_1 \to M_2$ is said to be \emph{Lipschitz of
order $\alpha > 0$} with constant $C \ge 0$ if
\begin{equation}
        d_2(f(x), f(y)) \le C \, d_1(x, y)^\alpha
\end{equation}
for every $x, y \in M_1$.  As before, this holds with $C = 0$ if and
only if $f$ is constant, and Lipschitz mappings of any order are
uniformly continuous.  If a real-valued function on the real line is
Lipschitz of order $\alpha > 1$, then it is constant, because it has
derivative $0$, although one could also show this more directly.  It
follows that a Lipschitz mapping of order $\alpha > 1$ from an
interval in the real line into any metric space is constant as well,
by composing with real-valued Lipschitz functions of order $1$ on the
range, such as the distance to a fixed point.

        Suppose that $0 < \beta < 1$.  If $r, t \ge 0$, then
\begin{equation}
        \max(r, t) \le (r^\beta + t^\beta)^{1/\beta}.
\end{equation}
Therefore
\begin{equation}
        r + t \le \max(r, t)^{1 - \beta} \, (r^\beta + t^\beta)
               \le (r^\beta + t^\beta)^{1/\beta},
\end{equation}
or equivalently
\begin{equation}
        (r + t)^\beta \le r^\beta + t^\beta.
\end{equation}
This is also very easy to check algebraically when $\beta = 1/2$, for
instance.

        If $(M, d(w, z))$ is a metric space, then it follows that
$d(w, z)^\beta$ is a metric on $M$ too when $0 < \beta < 1$.  This
does not work when $\beta > 1$, even for the real line.  Observe that
$f : M_1 \to M_2$ is Lipschitz of order $\alpha$ with respect to
$d_1(x, y)$ on $M_1$ if and only if $f$ is Lipschitz of order
$\alpha/\beta$ with respect to $d_1(x, y)^\beta$, keeping $d_2(u, v)$
fixed on $M_2$.  Similarly, $f$ is Lipschitz of order $\alpha$ with
respect to $d_2(u, v)$ on $M_2$ if and only if $f$ is Lipschitz of
order $\alpha \, \beta$ with respect to $d_2(u, v)^\beta$ on $M_2$,
keeping $d_1(x, y)$ fixed on $M_1$.

        A curve $p : [a, b] \to M$ in a metric space $(M, d(w, z))$
parameterized by a Lipschitz mapping of order $\alpha < 1$ can be
quite different from the case where $\alpha = 1$.  The length of $p$
can be infinite, and moreover $p([a, b])$ can be fractal.  This
includes common examples of snowflake curves in the plane.  Instead
one can show that the $\alpha$-dimensional Hausdorff measure of $p([a,
b])$ is finite.


\begin{thebibliography}{9}

\addcontentsline{toc}{section}{References}


\bibitem {b} R.~Beals, {\it Analysis: An Introduction}, Cambridge
University Press, 2004.

\bibitem {f} K.~Falconer, {\it The Geometry of Fractal Sets},
Cambridge University Press, 1986.

\bibitem {g} R.~Goldberg, {\it Methods of Real Analysis}, 2nd edition,
Wiley, 1976.

\bibitem {k} S.~Krantz, {\it Real Analysis and Foundations}, 2nd
edition, Chapman \& Hall / CRC, 2005.

\bibitem {pp} A.~Papadopoulos, {\it Metric Spaces, Convexity and
Nonpositive Curvature}, European Mathematical Society, 2005.

\bibitem {ros} M.~Rosenlicht, {\it Introduction to Analysis}, Dover,
1986.

\bibitem {rud} W.~Rudin, {\it Principles of Mathematical Analysis},
3rd edition, McGraw-Hill, 1976.

\bibitem {s1} S.~Semmes, {\it What is a metric space?},
arXiv:0709.1676 [math.MG].

\bibitem {s2} S.~Semmes, {\it An introduction to the geometry of
metric spaces}, arXiv:0709.4239 [math.MG].


\end{thebibliography}
\end{document}